\documentclass[12pt,reqno]{amsart}

\usepackage{amsmath}
\usepackage{amssymb}
\usepackage{amsthm}
\usepackage{graphicx}
\usepackage{tabularx}
\usepackage{array}
\usepackage{physics}
\usepackage{cleveref}
\crefname{thm}{Theorem}{Theorems}
\crefname{lem}{Lemma}{Lemmas}
\crefname{prop}{Proposition}{Propositions}
\crefname{cor}{Corollary}{Corollaries}

\DeclareMathOperator{\AFF}{aff}

\DeclareMathOperator{\EXT}{ext}

\DeclareMathOperator{\conv}{conv}

\DeclareMathOperator{\CH}{conv}
\DeclareMathOperator{\INT}{int}
\DeclareMathOperator{\BD}{bd}

\DeclareMathOperator{\RI}{relint}
\DeclareMathOperator{\REBD}{relbd}

\newcommand{\Espace}[1][n]{\mathbb{R}^{#1}}
\renewcommand{\norm}[1]{\left\Vert#1\right\Vert}
\newcommand{\origin}{o}

\newcommand{\SET}[2]{\left\{#1\mid #2\right\}}
\newcommand{\PI}{\mathbb{Z}^{+}}
\newcommand{\SCB}[1][n]{\mathcal{K}^{#1}}

\newcommand{\R}{\mathbb{R}}

\newtheorem{thm}{Theorem}
\newtheorem{cor}[thm]{Corollary}
\newtheorem{lem}[thm]{Lemma}

\newtheorem{conj}{Conjecture}

\theoremstyle{definition}

\newtheorem{rem}[thm]{Remark}

\begin{document}

\title[]{Homothetic covering of convex hulls of compact convex sets}


\author{Senlin Wu}

\email{wusenlin@nuc.edu.cn}

\author{Keke Zhang}

\email{S1908029@st.nuc.edu.cn}

\author{Chan He}

\email{hechan@nuc.edu.cn}

\address{Department of Mathematics, North University of China, 030051 Taiyuan,
  China}

\begin{abstract}
  Let $K$ be a compact convex set and $m$ be a positive integer. The covering
  functional of $K$ with respect to $m$ is the smallest $\lambda\in[0,1]$ such
  that $K$ can be covered by $m$ translates of $\lambda K$. Estimations of the
  covering functionals of convex hulls of two or more compact convex sets are
  presented. It is proved that, if a three-dimensional convex body $K$ is the
  convex hull of two compact convex sets having no interior points, then the
  least number $c(K)$ of smaller homothetic copies of $K$ needed to cover $K$ is
  not greater than $8$ and $c(K)=8$ if and only if $K$ is a parallelepiped.
\end{abstract}

\maketitle
\section{Introduction}
\label{sec:introduction}

Let $K$ be a compact convex subset of $\R^n$ that contains distinct points. We
denote by $\RI K$, $\REBD{K}$, $\INT K$, $\BD K$, and $\EXT{K}$ the
\emph{relative interior}, \emph{relative boundary}, \emph{interior},
\emph{boundary}, and the set of \emph{extreme points} of $K$, respectively. For
each $x\in\R^n$ and $\lambda\in(0,1)$, the set
\begin{displaymath}
  x+\lambda K:=\SET{x+\lambda y}{y\in K}
\end{displaymath}
is called a \emph{smaller homothetic copy} of $K$. We denote by $\SCB$ the set of
\emph{convex bodies} in $\Espace$, i.e., the set of compact convex sets in $\R^n$
having interior points. 

For each compact convex set $K$, we denote by $c(K)$ the least number of
translates of $\RI K$ needed to cover $K$. Concerning the least upper bound of
$c(K)$ in $\SCB$, there is a long standing conjecture (see \cite{Hadwiger1957},
\cite{Boltyanski-Martini-Soltan1997}, and \cite{Bezdek-Khan2018} for
more information about this conjecture):
\begin{conj}[Hadwiger's covering conjecture]\label{conj:main}
  For each $K\in\mathcal{K}^n$, $c(K)$ is bounded from the above by $2^n$, and
  this upper bound is attained only by parallelotopes.
\end{conj}
The assertion ``$c(K)\leq 2^n,~\forall K\in\SCB$'' will be referred to as the
``inequality part'' of Conjecture \ref{conj:main}. This conjecture has been
completely verified for several classes of convex bodies including: all planar
convex bodies (cf. \cite{Levi1955}), zonotopes, zonoids, belt bodies (cf.
\cite[\S 34]{Boltyanski-Martini-Soltan1997}), and convex hulls of a pair of
compact convex sets contained in two parallel hyperplanes in $\R^3$ (cf.
\cite{Wu-Zhou2019}). And the inequality part of Conjecture \ref{conj:main} has
been verified for centrally symmetric convex bodies in $\R^3$ (cf.
\cite{Lassak1984}), convex polyhedron in $\R^3$ having an affine symmetry (cf.
\cite{Bezdek1991}), convex bodies in $\R^3$ symmetric about a plane (cf.
\cite{Dekster2000}).

For each $m\in\PI$, we use the short-hand notation
\begin{displaymath}
  [m]=\SET{t\in\PI}{1\leq t\leq m}. 
\end{displaymath}
Note that, for each compact convex set $K$,
$c(K)$ equals the least number of smaller homothetic copies of $K$ needed to
cover $K$ (see, e.g., \cite[p. 262, Theorem 34.3]
{Boltyanski-Martini-Soltan1997}). Therefore, $c(K)\leq m$ for some $m\in\PI$ if
and only if $\Gamma_m(K)<1$, where $\Gamma_m(K)$ is defined by
\begin{displaymath}
  \Gamma_{m}(K):=\min\qty{\gamma>0\mid
  \exists\SET{x_{i}}{i\in[m]}\subseteq\Espace~\text{s.t.}~K\subseteq\bigcup\limits_{i=1}^m(x_{i}+\gamma K)},
\end{displaymath}
and is called \emph{the covering functional of $K$ with respect to $m$} (cf.
\cite{Lassak1986}, where $\Gamma_m(K)$ is called \emph{the $m$-covering number of $K$},
and \cite{Zong2010}). 

In this paper, we extend the results in \cite{Wu-Zhou2019} by studying the
homothetic covering problem for compact convex sets that can be expressed as
convex hulls of two or more compact convex sets. In Section 2, we provide an
estimation of covering functionals for this class of convex bodies in $\SCB$.
In Section 3, we solve Hadwiger's conjecture with respect to convex bodies in
$\SCB[3]$ that are convex hulls of two compact convex sets having empty
interiors.

\section{Covering functional of convex hulls of compact convex sets}
\label{sec:general-estimations}

The following estimation of the covering functionals of convex hulls of compact
convex sets $K_1,\dots,K_p$ uses only information about the covering functionals
of each $K_i$.

\begin{thm}
  \label{thm:first-estimation}
  Suppose that $K\in\SCB$ is the convex hull of convex compact sets
  $K_1,\dots,K_p$ and $m_1,\dots,m_p\in\PI$.
  \begin{enumerate}
  \item If $p\leq n+1$, then
    \begin{displaymath}
      \Gamma_{m_1+\dots+m_p}(K)\leq\max\SET{\frac{p-1+\Gamma_{m_i}(K_i)}{p}}{i\in[p]}.
    \end{displaymath}
  \item If $p>n+1$, then
    \begin{displaymath}
      \Gamma_{m_1+\dots+m_p}(K)\leq\max\SET{\frac{n+\Gamma_{m_i}(K_i)}{n+1}}{i\in[p]}.
    \end{displaymath}
  \end{enumerate}

\end{thm}
\begin{proof}
  Without loss of generality, we may assume that $\origin\in\RI{K}$. For each
  $i\in[p]$, put $\gamma_i=\Gamma_{m_i}(K_i)$. Then, for each $i\in [p]$, there
  exists a set $\SET{y_j^i}{j\in[m_i]}$ of $m_i$ points such that
  \begin{displaymath}
    K_i\subseteq \bigcup\limits_{j=1}^{m_i}(y_j^i+\gamma_i K_i)\subseteq \bigcup\limits_{j=1}^{m_i}(y_j^i+\gamma_i K).
  \end{displaymath}
  Let $x$ be an arbitrary point in $K$.

  {\bfseries Case 1}. $p\leq n+1$. By Theorem 3.13 in \cite{Soltan2015}, there
  exist $p$ points $x_1,\dots,x_p$, $p$ numbers
  $\lambda_1,\dots,\lambda_p\in[0,1]$ such that
  \begin{displaymath}
    x_i\in K_i,~\forall i\in[p],\quad
    \sum\limits_{i\in[p]}\lambda_i=1,\qqtext{and} x=\sum\limits_{i\in[p]}\lambda_ix_i.
  \end{displaymath}
  We may assume, without loss of generality, that
  \begin{displaymath}
    \lambda_1\geq \frac{1}{p}\qqtext{and} x_1\in y_1^1+\gamma_1 K_1\subseteq y_1^1+\gamma_1 K.
  \end{displaymath}
  Then
  \begin{align*}
    x=\lambda_1x_1+\sum\limits_{i=2}^{p}\lambda_ix_i&=\frac{1}{p}x_1+\qty(\lambda_1-\frac{1}{p})x_1+\sum\limits_{i=2}^{p}\lambda_ix_i\\
                                                    &\in\frac{1}{p}y_1^1+\frac{p-1+\gamma_1}{p}K\\
                                                    &\subseteq\frac{1}{p}y_1^1+\max\SET{\frac{p-1+\gamma_i}{p}}{i\in[p]}K.
  \end{align*}
  It follows that
  \begin{displaymath}
    K\subseteq \bigcup\limits_{i\in[p]}\SET{\frac{1}{p}y_j^i}{j\in[m_i]}+\max\SET{\frac{p-1+\gamma_i}{p}}{i\in[p]}K.
  \end{displaymath}

  {\bfseries Case 2}. $p>n+1$. By the Carath\'eodory's theorem, there exist
  $n+1$ points $x_1,\dots,x_{n+1}\in\bigcup\limits_{i=1}^pK_i$ and $n+1$ numbers
  $\lambda_1,\dots,\lambda_{n+1}\in[0,1]$ such that
  \begin{displaymath}
    x=\sum\limits_{i=1}^{n+1}\lambda_ix_i\qqtext{and} \sum\limits_{i=1}^{n+1}\lambda_i=1.
  \end{displaymath}
  We may assume, without loss of generality, that
  \begin{displaymath}
    \lambda_1\geq \frac{1}{n+1}\qqtext{and} x_1\in y_1^1+\gamma_1 K_1\subseteq y_1^1+\gamma_1 K.
  \end{displaymath}
  Then, in a similar way as above, we can show that
  \begin{displaymath}
    x\in\frac{1}{n+1}y_1^1+\max\SET{\frac{n+\gamma_i}{n+1}}{i\in[p]}K.
  \end{displaymath}
  It follows that
  \begin{displaymath}
    K\subseteq \bigcup\limits_{i\in[p]}\SET{\frac{1}{n+1}y_j^i}{j\in[m_i]}+\max\SET{\frac{n+\gamma_i}{n+1}}{i\in[p]}K.\qedhere
  \end{displaymath}
\end{proof}

In particular, we have the following:
\begin{cor}
  \label{cor:general-estimation}
  Suppose that $K$ is the convex hull of two non-empty compact convex sets $L$
  and $M$, and that $m_1,m_2\in\PI$. Then
  \begin{displaymath}
    \Gamma _{m_1+m_2}(K)\leq \max\qty{\frac{1+\Gamma_{m_1}(L)}{2},\frac{1+\Gamma_{m_2}(M)}{2}}.
  \end{displaymath}
\end{cor}

\begin{cor}
  \label{cor:segments}
  If $K$ is the convex hull of segments $K_1,\dots,K_p$, then
  \begin{displaymath}
    \Gamma_{2p}(K)\leq \frac{2p-1}{2p}.
  \end{displaymath}
\end{cor}

\begin{rem}
  When applying \cref{thm:first-estimation} to get a good estimation of
  $\Gamma_m(K)$, a suitable representation of $K$ as the convex hull of compact
  convex sets is necessary. For example, let $K$ be a three-dimensional simplex
  with $a,b,c,d$ as vertices. If we use the representation
  $K=\CH(\qty{a}\cup\CH\qty{b,c,d})$ then, by \cref{thm:first-estimation}, we
  have $\Gamma_4(K)\leq \frac{5}{6}$. But, if we use
  $K=\CH(\CH\qty{a,b}\cup\CH\qty{c,d})$, we will have the estimation
  $\Gamma_4(K)\leq \frac{3}{4}$, which is much better.

  When $n$ is odd and $K$ is an $n$-dimensional simplex, it is not difficult to
  verify that $\Gamma_{n+1}(K)=\frac{n}{n+1}$. By \cref{cor:segments}, we have
  $\Gamma_{n+1}(K)\leq \frac{n}{n+1}$. This shows that the estimation in
  \cref{thm:first-estimation} is tight in general. However, it can be improved
  in many other cases by taking the extremal structure of $K$ into consideration.
\end{rem}

\section{The three-dimensional case}

Let $K\subset\R^n$ be a compact convex set, $x\in\REBD K$, and $u\in\Espace$ be a
non-zero vector. If there exists a scalar $\lambda>0$ such that $x+\lambda u \in
\RI K$, then we say that $u$ \emph{illuminates} $x$. It is not difficult to see
that, a set $D$ of directions illuminates $\REBD K$ if and only if $D$ illuminates
all extreme points of $K$. Moreover (cf. Theorem 34.3 in
\cite{Boltyanski-Martini-Soltan1997}), $c(K)$ equals to the minimal cardinality
of a set of directions that can illuminate $\REBD K$.

A pair of points $a,b$ in a set $X\subseteq\R^n$ is called \emph{antipodal}
provided there are distinct parallel hyperplanes $H_a$ and $H_b$ through $a$ and
$b$, respectively, such that $X$ lies in the slab between $H_a$ and $H_b$.

\begin{lem}
  \label{lem:antipodal}
  Let $K\in\SCB$ and $x$, $y\in \BD K$. If $x$ and $y$ are not antipodal, then
  there is a direction that illuminates both $x$ and $y$.
\end{lem}
\begin{proof}
  We only need to consider the case when $x\neq y$. Since $x$ and $y$ are not
  antipodal, the segment $\qty[x,y]$ is not an affine diameter (cf.
  \cite{Soltan2005} for the definition and basic properties of affine diameters)
  of $K$. Let $\qty[u,v]$ be an affine diameter of $K$ parallel to $\qty[x,y]$
  and $c$ be an interior point of $K$. Without loss of generality, we may assume
  that
  \begin{displaymath}
    \frac{x-y}{\norm{x-y}}=\frac{u-v}{\norm{u-v}}.
  \end{displaymath}
  Then there exists a number $\lambda\in (0,1)$ such that $s-t=x-y$, where
  \begin{displaymath}
    s=\lambda c+(1-\lambda)u,\quad t=\lambda c+(1-\lambda)v.
  \end{displaymath}
  Clearly, both $s$ and $t$ are interior points of $K$. Let
  $d=\frac{s+t}{2}-\frac{x+y}{2}$. Then
  \begin{displaymath}
    x+d=\frac{x-y}{2}+\frac{s+t}{2}=s,\qqtext{and} y+d=\frac{y-x}{2}+\frac{s+t}{2}=t.
  \end{displaymath}
  I.e., $x$ and $y$ are both illuminated by $d\neq \origin$.
\end{proof}

\begin{thm}
  Let $K\in\SCB[3]$ be a convex body. If there exist two compact convex sets $L$
  and $M$ with empty interior such that $K=\conv(L\cup M)$, then $c(K)\leq 8$
  and the equality holds if and only if $K$ is a parallelepiped. 
\end{thm}
\begin{proof}
  We denote by $\AFF{L}$ and $\AFF{M}$ the affine dimensions of $L$ and $M$,
  respectively. We distinguish four cases.

  {\bfseries Case 1.} $0\in\qty{\AFF{L},\AFF{M}}$. Assume without loss of
  generality that $\AFF{L}=0$. Then $\AFF{M}=2$. By Theorem 4 in \cite
  {Wu-Xu2018} and the fact that $\Gamma_7(M)\leq\frac{1}{2}$ holds for each
  planar convex body (cf. \cite{Lassak1987}), we have
  \begin{displaymath}
    \Gamma _{8}(K) \leq \frac{1}{2-\Gamma_{7}(M)} \leq \frac{2}{3}.
  \end{displaymath}
  By \cref{cor:general-estimation} and the fact that $\Gamma_4(M)\leq
  \frac{\sqrt2}{2}$ holds for each planar convex body $M$ (cf.
  \cite{Lassak1986}), we have
  \begin{displaymath}
    \Gamma_5(K)\leq \frac{1+\frac{\sqrt{2}}{2}}{2}\approx 0.854.
  \end{displaymath}

  {\bfseries Case 2.} $\AFF{L}=\AFF{M}=1$. In this situation, $K$ is a
  three-dimensional simplex. We have (cf. \cite{Zong2010})
  \begin{displaymath}
    \Gamma_8(K)\leq \Gamma_5(K)=\frac{9}{13}.
  \end{displaymath}

  {\bfseries Case 3.} $\qty{\AFF{L},\AFF{M}}=\qty{1,2}$. Assume without loss of
  generality that $\AFF{L}=1$ and $\AFF{M}=2$. Then $\Gamma_2(L)=\frac{1}{2}$
  and $\Gamma_6(M)\leq \sin^{2}\frac{3\pi}{10}$ (cf. \cite{Lassak1987}). Then,
  \cref{cor:general-estimation} shows that
  \begin{displaymath}
    \Gamma_8(K)\leq \frac{1+\Gamma_{6}(M)}{2}\leq
    \frac{1+\sin^{2}\frac{3\pi}{10}}{2}\approx 0.827.
  \end{displaymath}
  In a similar way as in Case 1, we have
  \begin{displaymath}
    \Gamma_6(K)\leq \frac{1+\frac{\sqrt{2}}{2}}{2}\approx 0.854.
  \end{displaymath}
  
  {\bfseries Case 4. } $\AFF{L}=\AFF{M}=2$. In this case we have
  \begin{displaymath}
    \Gamma_8(K)\leq \frac{1+\frac{\sqrt{2}}{2}}{2}\approx 0.854.
  \end{displaymath}

  In the rest, we characterize the case when $c(K)=8$. The foregoing statements
  imply that $\AFF{L}=\AFF{M}=2$. If one of $L$ and $M$, say $L$, is not a
  parallelogram, then $\Gamma_3(L)<1$ and $\Gamma_4(M)<1$. It follows from
  \cref{cor:general-estimation} that
  \begin{displaymath}
    \Gamma_7(K)\leq \max\qty{\frac{1+\Gamma_3(L)}{2},\frac{1+\Gamma_4(M)}{2}}<1.
  \end{displaymath}
  Thus, $c(K)\leq 7$, a contradiction. In the following we assume that both $L$
  and $M$ are parallelograms. Since $c(K)=8$, $\EXT{K}=\EXT{L}\cup\EXT{M}$
  consisting of $8$ points. \cref{lem:antipodal} shows that the points in
  $\EXT{L}\cup\EXT{M}$ are pairwise antipodal. By the main result in
  \cite{Danzer-Grunbaum1962} (see also p. 225 in \cite{Martini-Soltan2005}),
  $\EXT{K}$ is the set of vertices of a parallelepiped. When $K$ is a
  parallelepiped, it is clear that $c(K)=8$.
\end{proof}


\def\polhk#1{\setbox0=\hbox{#1}{\ooalign{\hidewidth
  \lower1.5ex\hbox{`}\hidewidth\crcr\unhbox0}}} \def\cprime{$'$}
  \def\polhk#1{\setbox0=\hbox{#1}{\ooalign{\hidewidth
  \lower1.5ex\hbox{`}\hidewidth\crcr\unhbox0}}} \def\cprime{$'$}
\providecommand{\bysame}{\leavevmode\hbox to3em{\hrulefill}\thinspace}
\providecommand{\MR}{\relax\ifhmode\unskip\space\fi MR }
\providecommand{\MRhref}[2]{%
  \href{http://www.ams.org/mathscinet-getitem?mr=#1}{#2}
}
\providecommand{\href}[2]{#2}

\end{document}